\documentclass[11pt]{article}
\usepackage{color} 
\usepackage{ifpdf}
\ifpdf
    \usepackage[pdftex]{graphicx}
    \usepackage[pdftex]{hyperref}
    \hypersetup{
        unicode=false,          
        pdftoolbar=true,        
        pdfmenubar=true,        
        pdffitwindow=false,     
        pdfstartview={FitH},    
        pdftitle={AKG Article},      
        pdfauthor={Andrew Gillette},   
        pdfsubject={Mathematics},    
        pdfcreator={Andrew Gillette},  
        pdfproducer={Andrew Gillette}, 
        pdfkeywords={math, mathematics}, 
        pdfnewwindow=true,      
        colorlinks=true,        
        linkcolor=red,          
        citecolor=blue,         
        filecolor=magenta,      
        urlcolor=cyan           
    }

\else
    \usepackage{graphicx}
    \usepackage{pstricks}
    
    \newcommand{\href}[2]{#2}
\fi

\usepackage{amsmath,amsfonts,amssymb,algorithm,algorithmic,psfrag}
\usepackage{amsthm}
\usepackage{array}
\usepackage{psfrag}
\usepackage{subfigure}
\usepackage{bm}
\usepackage{color}
\usepackage[curve]{xypic}
\usepackage{bigstrut}

\newtheorem{theorem}{Theorem}

\newcommand{\be}{\textnormal{$\textbf{e}$}}

\newcommand{\bv}{\textnormal{$\textbf{v}$}}

\newcommand{\bx}{\textnormal{$\textbf{x}$}}
\newcommand{\hcurl}{H(\textnormal{curl})}
\newcommand{\hdiv}{H(\textnormal{div})}

\newcommand{\calP}{\mathcal P}
\newcommand{\calQ}{\mathcal Q}
\newcommand{\N} {\mathbb N}
\newcommand{\R} {\mathbb R}
\newcommand{\calR}{\mathcal R}
\newcommand{\calS}{\mathcal S}

\newcommand{\calU}{\mathcal U}
\newcommand{\calY}{\mathcal Y}

\newcommand{\sldeg}{\textnormal{sldeg}}
\newcommand{\tr}{\textnormal{tr}}

\newcommand{\raw}{\rightarrow}

\newcommand{\ds}{\displaystyle}

\newcommand{\p}{\partial}

\newcommand{\pyr}{\texttt{int}}
\newcommand{\bub}{\texttt{b}}

\newcommand{\ol}[1]{\overline{#1}}

\title{Serendipity and Tensor Product Affine Pyramid Finite Elements}

\author{Andrew Gillette\footnote{
Department of Mathematics,
University of Arizona,
Tucson, AZ, USA 85721. 
{\it agillette@math.arizona.edu}}}

\begin{document}
\maketitle

\begin{abstract}
\noindent
Using the language of finite element exterior calculus, we define two families of $H^1$-conforming finite element spaces over pyramids with a parallelogram base.
The first family has matching polynomial traces with tensor product elements on the base while the second has matching polynomial traces with serendipity elements on the base.
The second family is new to the literature and provides a robust approach for linking between Lagrange elements on tetrahedra and serendipity elements on affinely-mapped cubes while preserving continuity and approximation properties.
We define shape functions and degrees of freedom for each family and prove unisolvence and polynomial reproduction results.
\end{abstract}

\section{Introduction}

The pyramid geometry, known to all as one of the wonders of the ancient world, has proven to be an essential shape in the modern world of finite element modeling.
Three-dimensional geometries for studies of physical phenomena are typically built using meshes of either tetrahedral or hexahedral elements.
Tetrahedral elements allow great flexibility in representing intricate geometrical features, but the computational cost per element can become excessive if high order methods are required.
Hexahedral elements have easily exploitable computational advantages due to their tensor-product structure, however, this structure also limits their ability to mesh arbitrary geometries.
A best-of-both-worlds approach, which has been pursued with increasing interest in recent years~\cite{ADS2016,BRMHG2014,BCD2010,CWMRW2016,CWarb2015c,CWarb2015,CWarb2015b,CWarb2015a,FKDN2015,NP2012,NP2012a,WV2015}, uses hybrid meshes of tetrahedra, hexahedra, and pyramid geometries to balance computational efficiency with geometric flexibility.

In this paper, we use the language of finite element exterior calculus~\cite{AFW2006,AFW2010} to characterize two families of $H^1$-conforming finite element spaces over pyramids with a parallelogram base, one that is already known and a second that is new.
The first family, denoted $\calY^{-}_r\Lambda^0$, can be used to link $H^1$-conforming tensor product finite elements of order $r$ with $H^1$-conforming tetrahedral finite elements of order $r$ (i.e., with Lagrange elements).
The description and analysis of this family is greatly influenced by the work of Nigam and Philips \cite{NP2012,NP2012a}, which in turn builds on a great deal of prior mathematical and engineering work regarding pyramid finite elements.

\begin{figure}
\begin{center}
\sbox{\strutbox}{\rule{0pt}{0pt}}           
\begin{tabular}[.8\textwidth]{@{\extracolsep{\fill}} cccc}
\parbox{.22\textwidth}{\includegraphics[width=.22\textwidth]{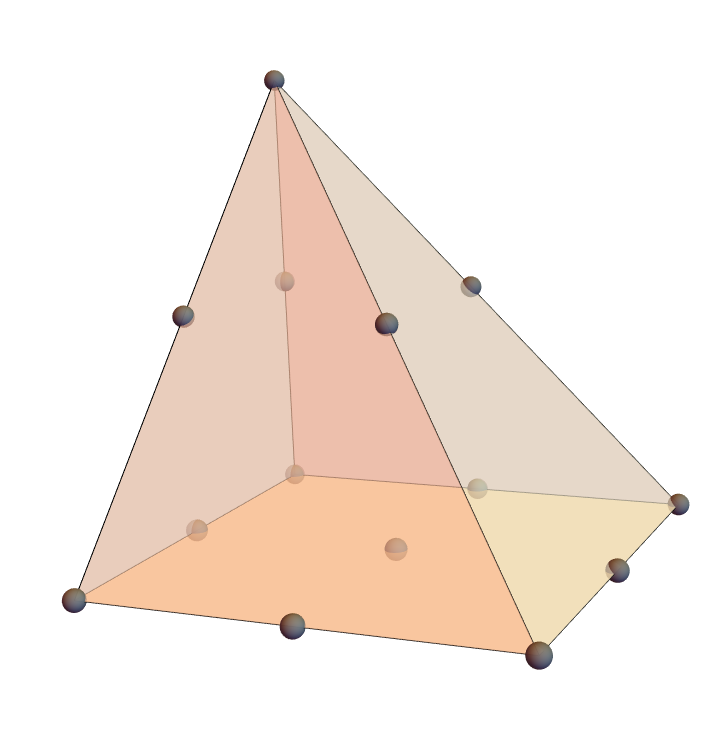}} & 
\parbox{.22\textwidth}{\includegraphics[width=.22\textwidth]{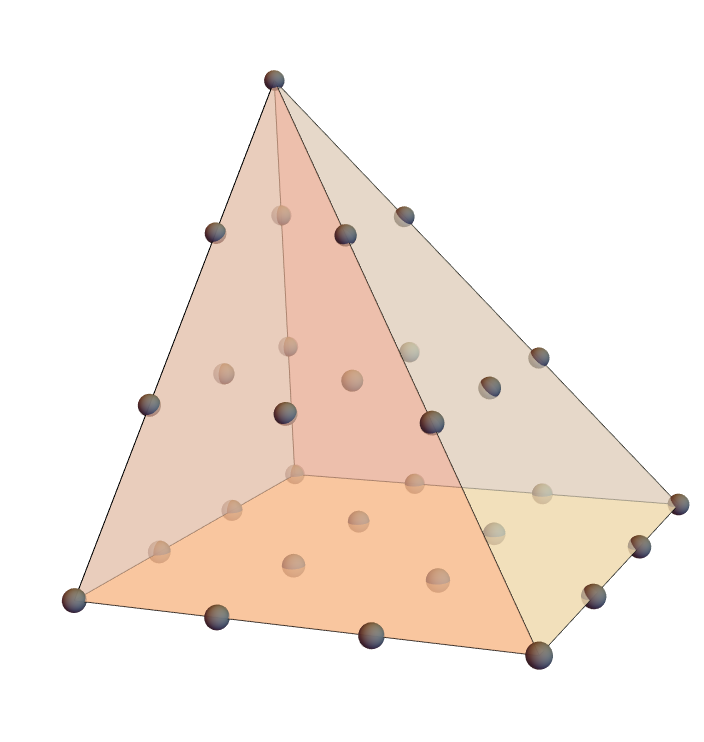}} &
\parbox{.22\textwidth}{\includegraphics[width=.22\textwidth]{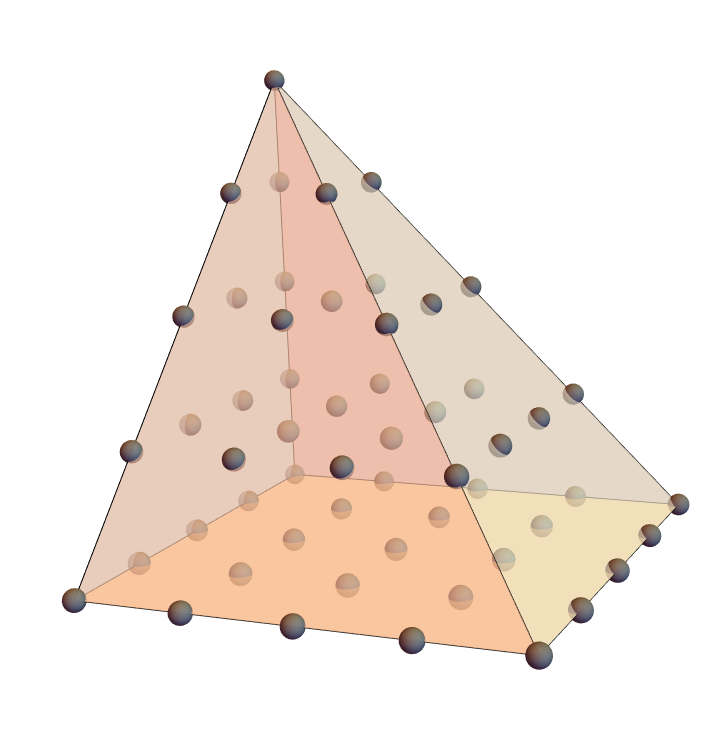}} & 
\parbox{.22\textwidth}{\includegraphics[width=.22\textwidth]{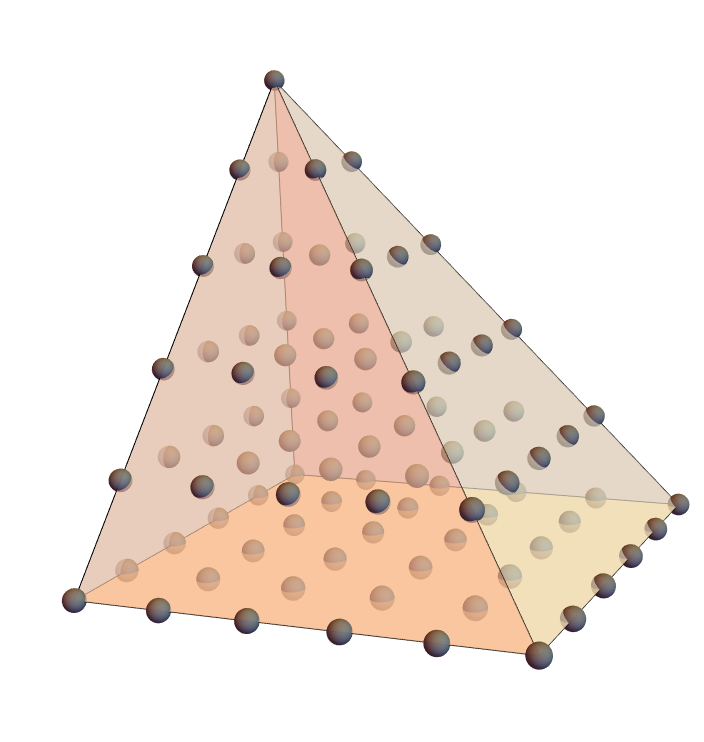}} \\
\parbox{.22\textwidth}{\includegraphics[width=.22\textwidth]{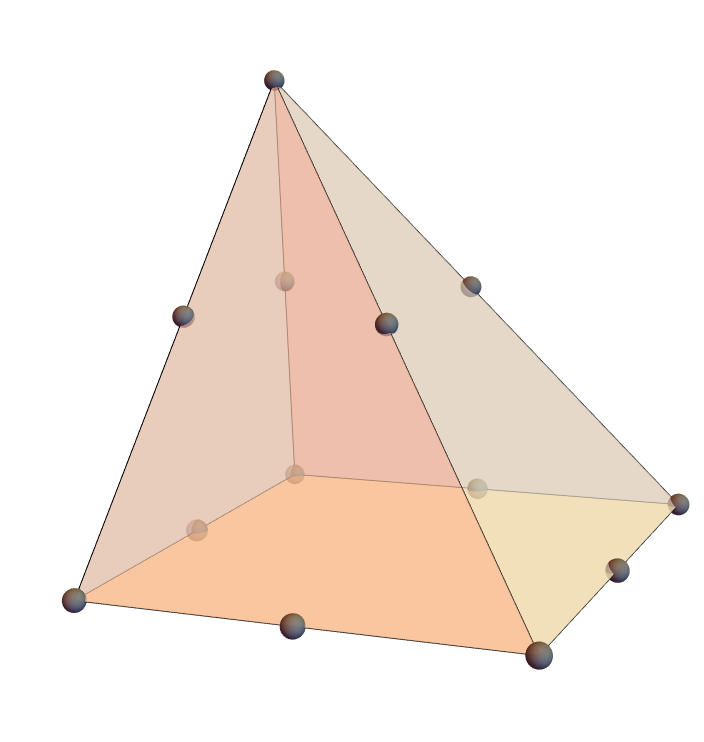}} & 
\parbox{.22\textwidth}{\includegraphics[width=.22\textwidth]{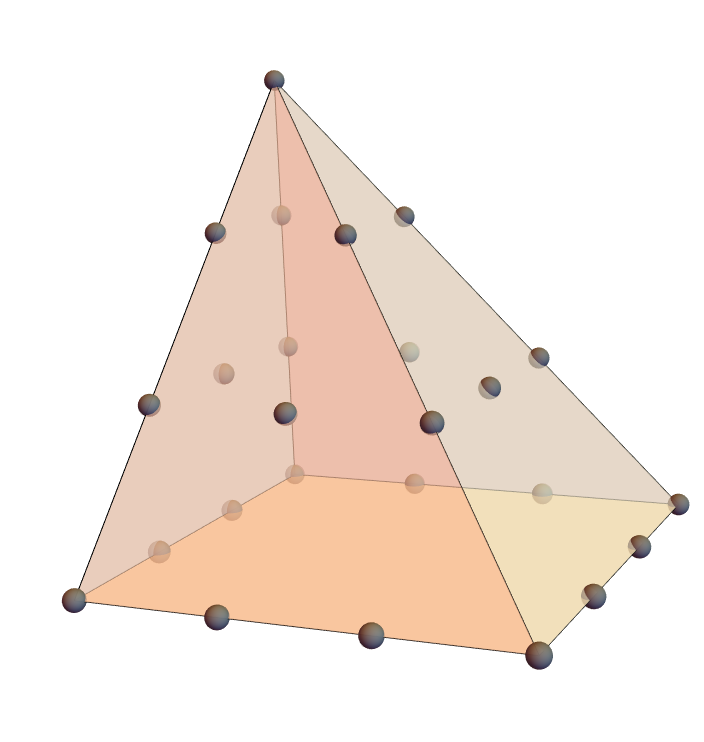}} &
\parbox{.22\textwidth}{\includegraphics[width=.22\textwidth]{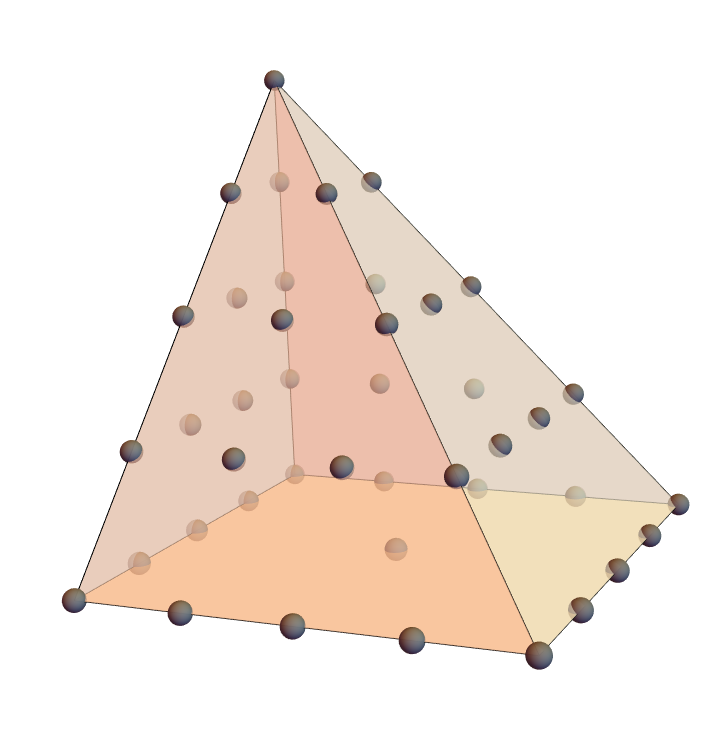}} & 
\parbox{.22\textwidth}{\includegraphics[width=.22\textwidth]{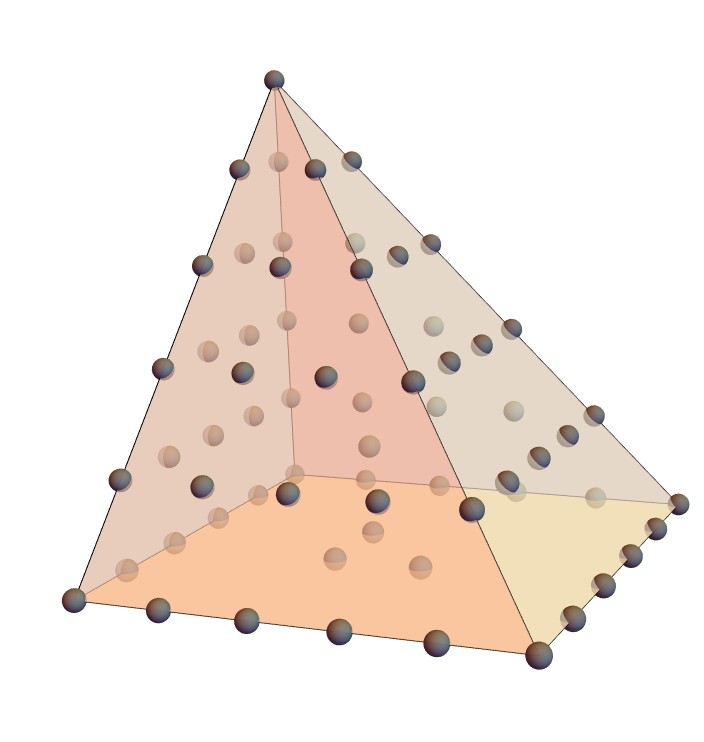}} \\
\end{tabular}
\end{center}
\caption{
Degrees of freedom for $\calY_r^-\Lambda^0$ (top row) and $\calY_r\Lambda^0$ (bottom row) for $r=$ 2, 3, 4, 5.
Each dot corresponds to a degree of freedom associated to the vertex, edge, face or interior where it is located.
See Section~\ref{sec:dim-opt} and Figure~\ref{fig:int-dofs} for more on the interior degrees of freedom.
}
\label{fig:dofs-viz}
\end{figure}

The second family, denoted $\calY_r\Lambda^0$, is new to the literature.
It can be used to link $H^1$-conforming serendipity finite elements of order $r$ on parallelpipeds with $H^1$-conforming tetrahedral finite elements of order $r$.
The definition makes use of the notion of the superlinear degree of a polynomial, as defined by Arnold and Awanou~\cite{AA2011}.
We show that
\[\dim \calY_r\Lambda^0 = \frac 16 (r^3+6r^2+23r)\leq \frac 16 (2 r^3+9 r^2+13 r+6) = \dim \calY^{-}_r\Lambda^0,\]
which becomes a strict inequality for any $r>1$.
Since serendipity elements provide a significant reduction in the number of degrees of freedom compared to tensor product elements, even for small $r$ values, the $\calY_r\Lambda^0$ family is poised to be a practical tool in the ongoing effort to minimize computational expense for problems on domains in $\R^3$.

Triangular prisms, sometimes called ``wedge'' elements, can also aid in bridging between tetrahedral and hexahedral meshes.
While helpful in specific meshing contexts, prisms have limited usefulness in a general context, due to the fact that their triangular faces occur on opposite sides of the geometry.
The faces of the pyramid are cleanly divided between the single quadrilateral base and the four triangular facets, making bridge-building between hexahedral and tetrahedral elements straightforward.
An additional consideration is the need to map physical mesh elements to reference elements via affine maps.
Any pyramid with a parallelogram base and apex not in the plane of the base can be mapped affinely back to a square-based reference pyramid.
Prisms, on the other hand, have three quadrilateral faces and thus require more constraints on when they can be mapped affinely to a reference element.

The outline of the remainder of this paper is as follows.
In Section~\ref{sec:bkgd}, we fix notation and explain prior work in this area.
Then, for each family, we present shape functions (Section~\ref{sec:shape-fns}), define degrees of freedom (Section~\ref{sec:dofs}), and prove results about unisolvence (Section~\ref{sec:unisolv}) and polynomial reproduction (Section~\ref{sec:polyn-rep}).
We conclude in Section~\ref{sec:dim-opt} with a discussion of dimension optimality.

\section{Background and Notation}
\label{sec:bkgd}

\paragraph{Finite element exterior calculus for simplices and cubes.}
Finite element exterior calculus~\cite{AFW2006,AFW2010} is a mathematical framework based on differential topology that describes and classifies many kinds of finite element methods in a unified fashion.
Appendix~\ref{appdx:feec} provides an introduction to some of the key ideas from finite element exterior calculus that are relevant to this work.
The Periodic Table of the Finite Elements~\cite{AL2014}, viewable online at \href{http://femtable.org}{femtable.org}, highlights many of the essential findings of the theory in visual form.

On a tetrahedron $\Delta_3$, the table identifies two families of elements, denoted $\calP_r^-\Lambda^k(\Delta_3)$ and $\calP_r\Lambda^k(\Delta_3)$.
In the scalar-valued case of interest here, the differential form order $k$ is $0$, and $\calP_r\Lambda^0(\Delta_3)$ is just the set of scalar-valued polynomials on $\Delta_3$ of degree at most $r$.
The spaces $\calP_r^-\Lambda^k(\Delta_3)$ are formally distinct from the $\calP_r\Lambda^k(\Delta_3)$ spaces.
However, in the $k=0$ and $k=3$ cases we have the identities $\calP_r^-\Lambda^0(\Delta_3)=\calP_r\Lambda^0(\Delta_3)$ and $\calP_r^-\Lambda^3(\Delta_3)=\calP_{r-1}\Lambda^3(\Delta_3)$.
We use both notations intentionally throughout the paper, as extensions of this work to higher $k$ values will require the distinction between these two families on simplices.

On a cube $\square_3$, the table identifies two families of elements, denoted $\calQ_r^-\Lambda^k(\square_3)$ and $\calS_r\Lambda^k(\square_3)$.
The $\calQ_r^-\Lambda^0(\square_3)$ family is the standard scalar-valued, tensor product element of order $r$, which has a total of $(r+1)^3$ degrees of freedom. 
The restriction of $\calQ_r^-\Lambda^0(\square_3)$ to a face $\square_2$ gives $\calQ_r^-\Lambda^0(\square_2)$, i.e.\ the tensor product element of order $r$, which has $(r+1)^2$ degrees of freedom.
The $\calS_r\Lambda^0(\square_3)$ family is the `serendipity' family of finite elements, known for some time in the mathematical and engineering literature~\cite{BS2002,Ci02,H1987,M1990,SF73,SB1991}, but generalized and characterized in a classical finite element setting only recently by Arnold and Awanou~\cite{AA2011}.
We present their definitions next.

\paragraph{Serendipity elements and superlinear degree.}

Given a multi-index $\alpha\in\N^n$, the \textit{degree} of the monomial $\bx^\alpha=\prod_{i=1}^n x_i^{\alpha_i}$ is $\deg(\bx^\alpha)=\sum_{i=1}^n\alpha_i$.
The \textit{superlinear degree} of $\bx^\alpha$, a term introduced in~\cite{AA2011}, is 
\[\sldeg (\bx^\alpha):= \sum_{\alpha_i\neq1}\alpha_i.\]
We note that $\sldeg(\bx^\alpha)\leq\deg(\bx^\alpha)$, with equality only when $\bx^\alpha$ has no variables that appear linearly.
The superlinear degree of a polynomial is the maximum of the superlinear degree of its monomials.
On an $n$-dimensional cube $\square_n$, the scalar-valued serendipity space is defined by
\begin{equation}
\label{eq:srdp-shpfns}
\calS_r\Lambda^0(\square_n) := \text{span}\left\{~\bx^\alpha~:~\sldeg(\bx^\alpha)\leq r~\right\}.
\end{equation}
In particular, note that for $r>0$, a basis for the the space $\calS_r\Lambda^0(\square_2)$ is given by the set of monomials $\{x^ay^b\}$ with $0\leq a+b\leq r$ or $(a,b)\in\{(1,r),(r,1)\}$.
The degrees of freedom for $\calS_r\Lambda^0(\square_n)$ are associated to its $d$-dimensional sub-faces $\square_d$ for $d=0,\ldots,n$, and are given by
\begin{equation}
\label{eq:srdp-dofs}
u\longmapsto\int_{\square_d} uq,\qquad \forall~q\in\calP_{r-2d}(\square_d),
\end{equation}
where $\calP_{r-2d}(\square_d)$ denotes polynomials on $\square_d$ of degree at most $r-2d$.
Arnold and Awanou proved in~\cite{AA2011} that the degrees of freedom from (\ref{eq:srdp-dofs}) are unisolvent for (\ref{eq:srdp-shpfns}).
Formally, (\ref{eq:srdp-dofs}) means
\begin{equation}
\label{eq:srdp-dofs-formal}
u\longmapsto\int_{\square_d} \left(\tr_{\square_d}u\right)\wedge q,\qquad \forall~q\in\calP_{r-2d}\Lambda^d(\square_d),
\end{equation}
where $\tr_{\square_d}u$ denotes the trace of $u$ on $\square_d$, and $\calP_{r-2d}\Lambda^d(\square_d)$ is the space of polynomial differential $d$-forms on $\square_d$ with coeffients in $\calP_{r-2d}$.
Additional background on traces and differential forms is given in Appendix~\ref{appdx:feec}.

\paragraph{Pyramid finite elements.}
The use of pyramid geometries in finite element methodologies began to gain attention with the work of Bedrosian~\cite{B1992}, Zgainski et al~\cite{ZCMCB1996}, and Coulomb et al~\cite{CZM1997} in the context of computational electromagnetics.
These and other early works focus primarily on questions related to implementation - an excellent summary is given in~\cite{BCD2010}.
More recently, Bergot, Cohen, and Durufl\'e~\cite{BCD2010} carried out a careful analysis of basis construction, interpolation error, and quadrature formulae for nodal pyramid elements of any polynomial order, including  physical elements that are non-affine maps of a reference element.
Nigam and Phillips~\cite{NP2012,NP2012a} allow only affinely-mapped reference elements, but provide compatibility, approximation, and stability results for $H^1$--, $\hcurl$--, and $\hdiv$--conforming pyramid elements in the context of exact sequences of finite element spaces.
Fuentes et al~\cite{FKDN2015} provide an implementation framework for the Nigam and Phillips elements as part of a complete $hp$-finite element package for ``hex-dominant'' meshes~\cite{BRMHG2014}.
Chan and Warburton have recently developed Berstein--B\'ezier style basis functions~\cite{CWarb2015b} and orthogonal bases~\cite{CWarb2015a} for the pyramid, as well as quadrature schemes~\cite{CWarb2015c}, trace inequalities~\cite{CWarb2015}, and implementations in discontinuous Galerkin settings (with additional collaborators)~\cite{CWMRW2016}.
A recent paper by by Ainsworth, Davydov and Schumaker~\cite{ADS2016} also looks at finite elements for tetrahedra-hexahedra-pyramid (THP) meshes with a view toward spline theory applications.

A classical finite element treatment of pyramid elements, meaning a triple of the form $\{$geometry, shape functions, degrees of freedom$\}$, has been rather elusive, the clearest examples appearing only recently in~\cite{BCD2010,NP2012,NP2012a}.
This is due in part to the fact that there are simple examples of low-degree polynomial traces on the faces of a pyramid that cannot be represented by a polynomial function satisfying the requisite inter-element compatibility criteria; a proof and discussion of this issue is given in the introduction of~\cite{NP2012}.
Proving that a space of rational functions is unisolvent for a set of degrees of freedom on the pyramid is sometimes handled indirectly by showing, for instance, that a Vandermonde matrix is invertible~\cite{BCD2010,CWarb2015c}.
A classical approach to proving unisolvency, given in ~\cite{NP2012} and~\cite{NP2012a}, presents degrees of freedom associated to the interior of a pyramid $\Omega$ as
\[u\longmapsto\int_\Omega\nabla u\cdot\nabla q~dV,\]
where $q$ belongs to the span of a set of rational bubble functions that are not explicitly stated.
Here, we present degrees of freedom associated to the interior of the pyramid that do not require a derivative on the input $u$ and generalize in a simple and obvious way to the new space $\calY_r\Lambda^0$.

Pyramid elements that link tetrahedral and serendipity elements for order $r>2$ have not been considered previously, to the best of our knowledge.
Liu et al~\cite{LDKG2011,LDYK2004} have defined sets of functions that might be used for the $r=2$ cases, however, the piecewise definition of these functions makes them computationally expensive and unlikely to generalize to $r>2$.

\paragraph{Reference geometries and mappings.}
We adopt the geometry conventions from Nigam and Phillips~\cite[Section 3.2]{NP2012a}, as restated here.
The \textit{infinite pyramid} element is
\[ K_\infty := \{~(x,y,z)\in\R^3\cup\infty~:~0\leq x,y\leq 1,~0\leq z\leq \infty~\}.\]
Formally, all points of the form $(x,y,\infty)$ are identified as a single point of $K_\infty$.
The \textit{reference pyramid} element is
\[ \hat K := \{~(\xi,\eta,\zeta)\in\R^3~:~0\leq \xi,\eta,\zeta,~ \xi\leq 1-\zeta,~ \eta\leq 1-\zeta~\}.\]
Thus, the four triangular faces of $\hat K$ are given by imposing one of the following additional constraints: $\xi=0$, $\eta=0$, $\xi=1-\zeta$, or $\eta=1-\zeta$.
The square base of $\hat K$ is given by imposing $\zeta=0$.
Define $\phi: K_\infty\raw \hat K$ by
\begin{align}
\phi(x,y,z) & = \begin{cases} \ds\left(\frac x{1+z},\frac y{1+z},\frac z{1+z}\right), & 0\leq z<\infty \\
(0,0,1), & z=\infty \end{cases} 
\end{align}
Note that $\phi$ suffices as a bijective change of coordinates between $(x,y,z)$ on $K_\infty$ and $(\xi,\eta,\zeta)$ on $\hat K$. 
In particular, we have that
\begin{equation}
\label{eq:phi-coord-chg}
\frac{x^ay^b}{(1+z)^c}~=~ \xi^a\eta^b(1-\zeta)^{c-a-b}.
\end{equation}
Given $u:\hat K\raw \R$, let $\phi^\ast u$ denote the pullback of $u$ to $K^\infty$ by $\phi$, that is:
\[\phi^\ast u:K_\infty\raw\R\quad\text{where}\quad (\phi^\ast u) (x,y,z) :=u(\phi(x,y,z)).\]
Since $\phi$ is a rational function in each coordinate, the pullback $\phi^\ast u$ of a polynomial function $u$ is, in general, a rational function.

An affine map of $\hat K$ will take the square base of $\hat K$ to a parallelogram embedded in $\R^3$.
Conversely, if $K\subset\R^3$ is a pyramid embedded in $\R^3$ with a parallelogram base, then there is an affine map that takes $\hat K$ to $K$.
The set of pyramids that are affine maps of $\hat K$ are called \textit{affine pyramids} by Nigam and Phillips~\cite{NP2012a} and we use the same terminology here.

\section{Shape Functions}
\label{sec:shape-fns}

\subsection{Shape functions for $\calY^-_r\Lambda^0$}

We now define spaces of rational shape functions on the infinite pyramid $K_\infty$ for the first family of pyramid finite elements, $\calY^-_r\Lambda^0$.
The construction here is exactly the same as Nigam and Phillips~\cite[Section 4.1]{NP2012a} and most of the notation is the same, with the notable difference that we use $r$ instead of $k$ to indicate polynomial degree.
Define 
\begin{equation}
\label{eq:Qrrr-def}
\calQ^{[r,r]}_r :=\text{span}\left\{\frac{x^ay^b}{(1+z)^c}~:~0\leq a,b\leq c\leq r~\right\}.
\end{equation}
We can decompose the space according to exponent of $(1+z)$ in the denominator.
This yields
\begin{equation}
\label{eq:Qrrr-decomp}
\calQ^{[r,r]}_r = \bigoplus_{j=0}^r \calQ^{j,j,0}_j \end{equation}
where
\begin{equation}
\label{eq:Qjj0j}
\calQ^{j,j,0}_j:=\text{span}\left\{\frac{x^ay^b}{(1+z)^j}~:~0\leq a,b\leq j~\right\}.
\end{equation}
Hence we have the dimension count:
\begin{equation}
\label{eq:qrrr-dim}
\dim \calQ^{[r,r]}_r = \sum_{j=0}^r\dim\calQ_j^-\Lambda^0(I^2) 
=\frac 16 (2 r^3+9 r^2+13 r+6).
\end{equation}
We define a set of shape functions on $\hat K$ by\rm
\footnote{A more precise notation for these spaces is $(\phi^{-1})^\ast\left(\calQ^{[r,r]}_r\right)$, as they are the set of pullbacks of $\calQ^{[r,r]}_r$ functions by $\phi^{-1}$; such notation is used by Nigam and Phillips.  We have used a simpler notation here only for the ease of reading.}
\begin{equation}
\label{eq:psh-fwd-qrrr}
\phi\left(\calQ^{[r,r]}_r\right) = \left\{~u:\hat K\raw~\R~:~\phi^\ast u\in \calQ^{[r,r]}_r~\right\}.
\end{equation}
Since $\phi$ is an isomorphism, $\dim\phi\left(\calQ^{[r,r]}_r\right) =\dim\calQ^{[r,r]}_r$.

\subsection{Shape functions for $\calY_r\Lambda^0$}

We now define spaces of rational shape functions on the infinite pyramid $K_\infty$ for the second family of pyramid finite elements, $\calY_r\Lambda^0$.
These spaces are \textit{new} to the literature, to the best of our knowledge, and fit naturally into the framework already developed.
The following definition was inspired from two key ideas: the decomposition of the shape function space $\calQ^{[r,r]}_r$ in terms of tensor product degrees of the numerator and the use of superlinear degree as a means of characterizing shape functions for serendipity spaces.
Define
\begin{equation}
\label{eq:Srrr-def}
\calS^{[r,r]}_r:= \text{span}\left\{\frac{x^ay^b}{(1+z)^c}~:~0\leq a,b\leq c \leq r,~~\sldeg(x^ay^b)\leq c ~\right\}.
\end{equation}
We can decompose the space according to exponent of $(1+z)$ in the denominator.
This yields
\begin{equation}
\label{eq:Srrr-decomp}
\calS^{[r,r]}_r = \bigoplus_{j=0}^r \calS^{j,j,0}_j
\end{equation}
where
\begin{equation}
\label{eq:Sjj0j}
\calS^{j,j,0}_j:=\text{span}\left\{\frac{x^ay^b}{(1+z)^j}~:~0\leq a,b\leq j,~~\sldeg(x^ay^b)\leq j~\right\}.
\end{equation}
Observe that for $j>0$, the constraint that $\sldeg(x^ay^b)\leq j$ implies that either $\deg(x^ay^b)\leq j$ or $x^ay^b\in\{xy^j,x^jy\}$.
Thus, if $j\not= 0$, we have
\[\dim\calS^{j,j,0}_j:=\dim\calS_j\Lambda^0(I^2)=\dim \calP_j\Lambda^0(I^2)+2 = \binom{2+j}{j} + 2,\]
which follows from, e.g.~\cite[Equation (2.1)]{AA2011}.
Since $\dim\calS^{0,0,0}_0=1$, we have the dimension count:
\begin{equation}
\label{eq:srrr-dim}
\dim \calS^{[r,r]}_r = 1 + \sum_{j=1}^r\dim \calS_j\Lambda^0(I^2) =\frac 16 (r^3+6r^2+23r).
\end{equation}
We define a set of shape functions on $\hat K$ by
\begin{equation}
\label{eq:psh-fwd-srrr}
\phi\left(\calS^{[r,r]}_r\right) = \left\{~u:\hat K\raw~\R~:~\phi^\ast u\in \calS^{[r,r]}_r~\right\}.
\end{equation}
Since $\phi$ is an isomorphism, $\dim\phi\left(\calS^{[r,r]}_r\right) =\dim\calS^{[r,r]}_r$.

\subsection{The ``lowest order bubble function'' on $K_\infty$ and $\hat K$}

A key function of interest to our subsequent analysis is $\bub:K_\infty\raw\R$ given by
\begin{equation}
\label{eq:bub-def}
\bub(x,y,z):=\frac{x(1-x)y(1-y)z}{(1+z)^3}.
\end{equation}
The numerator of $\bub$ indicates that it vanishes on the five `faces' of $K_\infty$ while the denominator indicates that $\bub(x,y,z)\raw 0$ as $z\raw\infty$.
Hence, $\bub$ vanishes on $\p K_\infty$. 
We can write
\[\bub(x,y,z)=\frac{x(1-x)y(1-y)}{(1+z)^{2}}-\frac{x(1-x)y(1-y)}{(1+z)^3}\in\calQ^{[3,3]}_3.\]
We call $\bub$ the ``lowest order bubble function'' on $K_\infty$ as the space $\calQ^{[r,r]}_r$ does not contain any functions that vanish identically on $\p K_\infty$ if $r<3$.
Changing coordinates by $\phi$  gives
\[\bub(\xi,\eta,\zeta) = \frac{\xi \eta \zeta (\xi+\zeta-1) (\eta+\zeta-1)}{(\zeta-1)^2}\in\phi\left(\calQ^{[3,3]}_3\right).\]
Note also that $\bub\in\calS^{[5,5]}_5$ and that the space $\calS^{[r,r]}_r$ does not contain any functions that vanish identically on $\p K_\infty$ if $r<5$.

\section{Degrees of Freedom}
\label{sec:dofs}

We now state all the degrees of freedom precisely in a classical finite element sense and count them.
Degrees of freedom for a function $u:\hat K\raw\R$ are defined in terms of its trace on vertices, edges, triangular faces, parallelogram face, and the interior of $\hat K$, integrated against functions from index spaces denoted $P_\bv$, $P_\be$, $P_\triangle$, $P_{\,\square}$,  $R_{\,\pyr}$, respectively.
The index spaces are defined in Table~\ref{tab:all-dofs-0}.
Since we are describing degrees of freedom for spaces of $0$-forms, the index space associated to a $d$-dimensional object is a space of differential $d$-forms; this point is discussed further in Appendix~\ref{appdx:feec}.
The spaces $R_{\,\pyr}$ are differential 3-forms with rational functions as coefficients while all the other spaces have polynomial coefficients.

To each vertex $\bv$ of the pyramid, associate the evaluation degree of freedom
\begin{equation}
\label{eq:dof-vtx}
u\longmapsto u(\bv).
\end{equation}
To each edge $\be$ of the pyramid, associate
\begin{equation}
\label{eq:dof-edge}
u\longmapsto \int_{\be}(\tr_\be\,u)\ q, \quad q\in P_\be.
\end{equation}
To each triangular face $\triangle$ of the pyramid, associate 
\begin{equation}
\label{eq:dof-tri}
u\longmapsto \int_{\triangle}(\tr_\triangle\,u)\ q, \quad q\in P_\triangle.
\end{equation}
To the parallelogram face $\square$ of the pyramid, associate
\begin{equation}
\label{eq:dof-quad}
u\longmapsto \int_{\square}(\tr_\square\,u)\ q, \quad q\in P_{\,\square}.
\end{equation}
To the three-dimensional interior $\pyr$ of the pyramid, associate
\begin{equation}
\label{eq:dof-int}
u\longmapsto \int_{\pyr}(\tr_\pyr\,u)\ q, \quad q\in R_{\,\pyr}.
\end{equation}

\begin{table}[h]
\begin{center}
$
\begin{array}{r|ccccc}
 & P_\bv & P_\be & P_\triangle & P_{\,\square} & R_{\,\pyr} \\
\hline 
&\\[-3mm]
\calY^{-}_r\Lambda^0 & \R & \calP_{r-2}\Lambda^{1}(\be) & \calP_{r-3}\Lambda^{2} (\triangle) & \calQ_{r-1}^-\Lambda^{2}(\square) &  \phi\left(\bub\cdot\calQ^{[r-3,r-3]}_{r-3}\right)\Lambda^3(\pyr)  \\[2mm]
\calY^{}_r\Lambda^0 & \R & \calP^-_{r-1}\Lambda^{1}(\be) & \calP^-_{r-2}\Lambda^{2}(\triangle) & \calP_{r-3}^-\Lambda^{2}(\square) &  \phi\left(\bub\cdot\calS^{[r-5,r-5]}_{r-5}\right)\Lambda^3(\pyr)
\end{array}
$
\end{center}
\caption{Index spaces for the degrees of freedom for the two pyramid families.
}
\label{tab:all-dofs-0}
\end{table}

Note that $\calP^-_{r}\Lambda^n(\R^n)=\calP_{r-1}\Lambda^{n}(\R^n)$, meaning the spaces for $P_\bv$, $P_\be$ and $P_\triangle$ are the same for both families.
On the parallelogram face $\square_2$, we recognize $\calQ_{r-1}^-\Lambda^{2}(\square_2)$ as the  indexing space for $\calQ^-_r\Lambda^0(\square_2)$, as expected for the $\calY_r^-\Lambda^0$ family.
Likewise, for the parallelogram face in the $\calY_r\Lambda^0$ family, we apply the identity  $\calP_{r-3}^-\Lambda^{2}(\square_2)=\calP_{r-4}\Lambda^2(\square_2)$ and recover the index space for $\calS_r\Lambda^0(\square_2)$.
For the spaces $R_\pyr$, the notation means
\[\phi\left(\bub\cdot\calQ^{[r-3,r-3]}_{r-3}\right)\Lambda^3(\pyr):= \text{span}\left\{~u\, dV~:~\phi^\ast u=\bub q~\text{with}~q\in \calQ^{[r-3,r-3]}_{r-3}~\right\},\]
where $dV$ is the volume 3-form $d\xi d\eta d\zeta$ on $\hat K$.
The meaning of $R_\pyr$ in the $\calY_r\Lambda^0$ case is analogous.

Table~\ref{tab:all-dofs-0} can be used to compute $\dim\calY^{-}_r\Lambda^0$ and $\dim\calY^{}_r\Lambda^0$ as follows.
We compute the dimension of each entry in Table~\ref{tab:all-dofs-0}, weight by the number of times each kind of geometrical object appears in the pyramid, and then sum.
This gives 
%
%
\begin{align}
\dim\calY^{-}_r\Lambda^0 & = 5+ 8 \left|\calP_{r-2}\Lambda^{1}(\be)\right| + 4\left|\calP_{r-3}\Lambda^{2}(\triangle)\right| + \left|\calQ_{r-1}^-\Lambda^{2}(\square)\right| + \left| \calQ^{[r-3,r-3]}_{r-3}\right| \notag\\
& = 5 + 8(r-1) + 2 (r-2) (r-1) + (r-1)^2 + \frac{(2r-3)(r-2)(r-1)}{6}\notag \\
& = \frac 16 (2 r^3+9 r^2+13 r+6). \label{eq:yrm-dim}
\end{align}

\begin{align}
\dim\calY_r\Lambda^0 & = 5+ 8 \left|\calP^-_{r-1}\Lambda^{1}(\be)\right| + 4\left|\calP^-_{r-2}\Lambda^{2}(\triangle)\right| + \left|\calP_{r-4}\Lambda^{2}(\square)\right| + \left|\calS^{[r-5,r-5]}_{r-5}\right|\notag \\
& = 5 + 8(r-1) + 2 (r-2) (r-1) + \frac{(r-3)(r-2)}{2} + 
\frac{(r-4) (r-3) (r-2)}{6}\notag \\
& = \frac 16 (r^3+6r^2+23r). \label{eq:yr-dim}
\end{align}
Comparing (\ref{eq:yrm-dim}) to (\ref{eq:qrrr-dim}) and (\ref{eq:yr-dim}) to (\ref{eq:srrr-dim}), we see that
\begin{equation}
\label{eq:dim-match}
\dim\calY_r^-\Lambda^0 = \dim  \calQ^{[r,r]}_r \quad\text{and}\quad \dim\calY_r\Lambda^0 = \dim  \calS^{[r,r]}_r.
\end{equation}

The degrees of freedom are associated to portions of the pyramid geometry according to their definition in (\ref{eq:dof-vtx})-(\ref{eq:dof-int}).
We visualize this association for $\calY_r^-\Lambda^0$ and $\calY_r\Lambda^0$ for $r=$ 2, 3, 4, 5 in Figure~\ref{fig:dofs-viz}.

\section{Unisolvence and $H^1$-conformity}
\label{sec:unisolv}

We now prove that the degrees of freedom for $\calY^{-}_r\Lambda^0$ and $\calY_r\Lambda^0$ are unisolvent for $\phi\left(\calQ^{[r,r]}_r\right)$ and $\phi\left(\calS^{[r,r]}_r\right)$, respectively.
As part of the proof, we show that given $u:\hat K\raw\R$ in one of the shape function spaces, the trace of $u$ on each boundary face of the pyramid is a bivariate polynomial.
Moreover, these polynomials have total degree at most $r$ on triangular facets, degree at most $r$ in each variable on the parallelogram face for $\calY_r^-\Lambda^0$, and superlinear degree at most $r$ on the parallelogram face for $\calY_r\Lambda^0$.
As a consequence, both $\calY^{-}_r\Lambda^0$ and $\calY_r\Lambda^0$ are guaranteed to be $H^1$-conforming when linked with tetrahedral and hexahedral elements of the corresponding types.

\newpage
\begin{theorem}[Unisolvence]
By employing the definitions from Table~\ref{tab:all-dofs-0}:
\begin{enumerate} 
\renewcommand{\theenumi}{\roman{enumi}}
\item The degrees of freedom for $\calY^{-}_r\Lambda^0$ are unisolvent for $\phi\left(\calQ^{[r,r]}_r\right)$;
\item The degrees of freedom for $\calY_r\Lambda^0$ are unisolvent for $\phi\left(\calS^{[r,r]}_r\right)$.
\end{enumerate}
\end{theorem}

\begin{proof}
We prove \textit{i.} first.
Since $\phi$ is a bijection, it follows from (\ref{eq:dim-match}) that $\dim\calY_r^-\Lambda^0 = \dim \phi\left( \calQ^{[r,r]}_r\right)$.
Let $u\in\phi\left( \calQ^{[r,r]}_r\right)$ and suppose that all the quantities in (\ref{eq:dof-vtx})-(\ref{eq:dof-int}) vanish, using the definitions from the top row of Table~\ref{tab:all-dofs-0}.
It suffices to show that $u$ vanishes.
Using (\ref{eq:Qrrr-def}) and (\ref{eq:phi-coord-chg}), we have that
\[u\in \text{span}\left\{\xi^a\eta^b(1-\zeta)^{c-a-b}~:~0\leq a,b\leq c\leq r~\right\}. \]
First we show that $u$ vanishes on each face.
On the quadrilateral face $\square$, we have $\zeta=0$, so that
\begin{equation}
\label{eq:uni-Q-sq}
\tr_\square u\in \text{span}\left\{\xi^a\eta^b~:~0\leq a,b\leq r~\right\} = \calQ_r^-\Lambda^0(\square).
\end{equation}
The degrees of freedom associated to $\square$ and its edges and vertices are unisolvent for $\calQ_r^-\Lambda^0(\square)$, so $u$ vanishes on $\square$.
For the triangular face with $\xi=0$, call it $\triangle_1$, we have
\begin{equation}
\label{eq:uni-Q-tri1}
\tr_{\triangle_1} u\in \text{span}\left\{\eta^b(1-\zeta)^{c-b}~:~0\leq b\leq c\leq r~\right\} = \calP_r^-\Lambda^0(\triangle_1).
\end{equation}
The degrees of freedom associated to $\triangle_1$ and its edges and vertices are unisolvent for $\calP_r^-\Lambda^0(\triangle_1)$, so $u$ vanishes on $\triangle_1$.
The other triangular faces follow similarly.  For instance, on the triangular face with $\xi=1-\zeta$, call it $\triangle_2$, we have
\begin{equation}
\label{eq:uni-Q-tri2}
\tr_{\triangle_2} u\in \text{span}\left\{\xi^{c-b}\eta^b~:~0\leq b\leq c\leq r~\right\} = \calP_r^-\Lambda^0(\triangle_2).
\end{equation}
Thus, it remains to show that if $u$ vanishes on $\square$ and on $\triangle_1-\triangle_4$ then $u\equiv 0$.
Now, since $u$ vanishes on $\p\hat K$, we have a function $\phi^\ast u\in \calQ^{[r,r]}_r$ that vanishes on $\p K_\infty$.
Write
\begin{equation}
\label{eq:phiu-expr-1}
\phi^\ast u=\sum_{i=0}^{r}\frac{p_i(x,y)}{(1+z)^i}
\end{equation}
for some $p_i$, polynomials in $x$ and $y$, with $p_i\in\calQ_{i}^-\Lambda^0(\R^2)$.
Since $\phi^\ast u$ is a polynomial in $x$ and $y$, that vanishes on $\{x=0\}$, $\{x=1\}$, $\{y=0\}$, and $\{y=1\}$, we can factor $x(1-x)y(1-y)$ out of the expression (\ref{eq:phiu-expr-1}).
This forces $p_0=p_1=0$ in (\ref{eq:phiu-expr-1}) as $x(1-x)y(1-y)\in\calQ_2^-\Lambda^0(\R^2)$ will not factor out of any function in $\calQ^{[0,0]}_0$ or $\calQ^{[1,1]}_1$.
Thus, we can write
\begin{equation}
\label{eq:phiu-expr-2}
\phi^\ast u=\frac{x(1-x)y(1-y)}{(1+z)^2}\sum_{i=0}^{r-2}\frac{\ol p_i(x,y)}{(1+z)^i}
\end{equation}
for some $\ol p_i$, polynomials in $x$ and $y$, with $\ol p_i\in\calQ_{i}^-\Lambda^0(\R^2)$.
Observe that no non-zero element of $\calQ^{[0,0]}_0$ vanishes on $\{z=0\}$, however
\begin{equation}
\label{eq:z-factor}
\frac{z}{1+z}=1-\frac{1}{1+z}\in\calQ^{[1,1]}_1
\end{equation}
does vanish on $\{z=0\}$.
Further, an element of  $\calQ^{[r,r]}_r$ will vanish on $\{z=0\}$ if and only if it is divisible by $\frac{z}{1+z}$, meaning we can write
\begin{equation}
\label{eq:phiu-expr-3}
\phi^\ast u=\frac{x(1-x)y(1-y)z}{(1+z)^3}\sum_{i=0}^{r-3}\frac{\hat p_i(x,y)}{(1+z)^i}
\end{equation}
for some $\hat p_i$, polynomials in $x$ and $y$, with $\hat p_i\in\calQ_{i}^-\Lambda^0(\R^2)$.
Hence, $u\in \phi\left(\bub\cdot\calQ^{[r-3,r-3]}_{r-3}\right)$ and we may take $q=u dV$ in (\ref{eq:dof-int}) to get
\[\int_{\hat K} u^2 dV=0.\]
Thus, $u=0$.

The proof of \textit{ii.}\ is similar.
Since $\phi$ is a bijection, it follows from (\ref{eq:dim-match}) that $\dim\calY_r\Lambda^0 = \dim \phi\left( \calS^{[r,r]}_r\right)$.
Let $u\in\phi\left( \calS^{[r,r]}_r\right)$ and suppose that all the quantities in (\ref{eq:dof-vtx})-(\ref{eq:dof-int}) vanish, using the definitions from the bottom row of Table~\ref{tab:all-dofs-0}.
Using (\ref{eq:Srrr-def}) and (\ref{eq:phi-coord-chg}), we have that
\begin{equation}
\label{eq:u-srrr-span}
u\in \text{span}\left\{\xi^a\eta^b(1-\zeta)^{c-a-b}~:~0\leq a,b\leq c\leq r, ~~\sldeg(\xi^a \eta^b)\leq c~\right\}.
\end{equation}
First we show that $u$ vanishes on each face.
On the quadrilateral face $\square$, we have $\zeta=0$, so that
\begin{equation}
\label{eq:uni-S-sq}
\tr_\square u\in \text{span}\left\{\xi^a\eta^b~:~\sldeg(\xi^a \eta^b)\leq r~\right\} = \calS_r\Lambda^0(\square).
\end{equation}
Note that since $r\geq 1$, the constraint $\sldeg(\xi^a \eta^b)\leq r$ ensures that $0\leq a,b\leq r$.
The degrees of freedom associated to $\square$ and its edges and vertices are unisolvent for $\calS_r^-\Lambda^0(\square)$ (see~\cite{AA2011} for a proof) so $u$ vanishes on $\square$.
For the triangular face with $\xi=0$, call it $\triangle_1$, take $a=0$ in (\ref{eq:u-srrr-span}), giving
\footnote{Recall that $\calP_r\Lambda^0=\calP_r^-\Lambda^0$; see the beginning of Section~\ref{sec:bkgd} for a comment on this.}
\begin{equation}
\label{eq:uni-S-tri1}
\tr_{\triangle_1} u\in \text{span}\left\{\eta^b(1-\zeta)^{c-b}~:~0\leq b\leq c\leq r~\right\} = \calP_r\Lambda^0(\triangle_1).
\end{equation}
The degrees of freedom associated to $\triangle_1$ and its edges and vertices are unisolvent for $\calP_r\Lambda^0(\triangle_1)$, so $u$ vanishes on $\triangle_1$.
The other triangular faces follow similarly. 
Thus, it remains to show that if $u$ vanishes on $\square$ and on $\triangle_1-\triangle_4$ then $u\equiv 0$.
Now, since $u$ vanishes on $\p\hat K$, we have a function $\phi^\ast u\in \calS^{[r,r]}_r$ that vanishes on $\p K_\infty$.  
As in the proof of part \textit{i.}, $x(1-x)y(1-y)$ must factor out of $\phi^\ast u$, and we note that $r=4$ is the smallest value of $r$ for which $x(1-x)y(1-y)\in\calS_r\Lambda^0(\R^2)$.
Thus, we can write
\begin{equation}
\label{eq:phiu-expr-s2}
\phi^\ast u=\frac{x(1-x)y(1-y)}{(1+z)^4}\sum_{i=0}^{r-4}\frac{\ol p_i(x,y)}{(1+z)^i}
\end{equation}
for some $\ol p_i$, polynomials in $x$ and $y$, with $\ol p_i\in\calS_{i}^-\Lambda^0(\R^2)$, interpreting $\calS_{0}^-\Lambda^0(\R^2)$ as $\R$.
Recalling (\ref{eq:z-factor}), an element of $\calS^{[r,r]}_r$ will vanish on $\{z=0\}$ if and only if it is divisible by $\frac z{1+z}$.
Thus,
\begin{equation}
\label{eq:phiu-expr-s3}
\phi^\ast u=\frac{x(1-x)y(1-y)z}{(1+z)^5}\sum_{i=0}^{r-5}\frac{\hat p_i(x,y)}{(1+z)^i}
\end{equation}
for some $\hat p_i$, polynomials in $x$ and $y$, with $\hat p_i\in\calS_{i}^-\Lambda^0(\R^2)$.
Hence, $u\in \phi\left(\bub\cdot\calS^{[r-5,r-5]}_{r-5}\right)$.
Therefore, we may again take $q=u dV$ in (\ref{eq:dof-int}) so that $u^2dV$ has integral 0 over $\hat K$ and thus $u=0$.
\end{proof}

\section{Polynomial reproduction and error analysis}
\label{sec:polyn-rep}

\begin{theorem}[Polynomial reproduction]
\label{thm:polyn-rep}
Let $\calP_r(\R^3)$ denote polynomials of degree at most $r$ in three variables.  We have:
\begin{enumerate} 
\renewcommand{\theenumi}{\roman{enumi}}
\item $\calP_r(\R^3)\subset\phi\left(\calQ^{[r,r]}_r\right)$,
\item $\calP_r(\R^3)\subset\phi\left(\calS^{[r,r]}_r\right)$.
\end{enumerate}
\end{theorem}

\begin{proof}
Let $p(\xi,\eta,\zeta)\in\calP_r(\R^3)$.
Expand $p$ in powers of $\xi$, $\eta$, and $(1-\zeta)$; this does not change the degree of $p$ nor its degree with respect to any variable. 
Without loss of generality, assume that this expansion yields a single monomial, i.e.\ $p=\xi^a \eta^b (1-\zeta)^{c}$ with $0\leq a,b\leq a+b+c\leq r$. 
From (\ref{eq:phi-coord-chg}), 
\[\phi^\ast p = \frac{x^ay^b}{(1+z)^{a+b+c}}\in\calQ^{[r,r]}_r.\]
We have established that $p$ is a real-valued function on $\hat K$ whose pullback, $\phi^\ast p$, is in $\calQ^{[r,r]}_r$.
By the definition of $\phi\left(\calQ^{[r,r]}_r\right)$ in (\ref{eq:psh-fwd-qrrr}), we have $p\in\phi\left(\calQ^{[r,r]}_r\right)$.
Further, $\sldeg(\xi^a\eta^b)\leq\deg(\xi^a\eta^b)=a+b\leq a+b+c\leq r$, so $p\in\phi\left(\calS^{[r,r]}_r\right)$.
\end{proof}

In regards to \textit{a priori} error estimates, it was shown by Arnold, Boffi, and Bonizzoni~\cite{ABB2012} that $\calS_r\Lambda^0$ elements on $n$-dimensional cubes have a standard $O(h^r)$ convergence rate when all physical mesh elements are affine maps of the reference element.
Similarly, Nigam and Phillips~\cite{NP2012a} and Bergot, Cohen, and Durufl\'e~\cite{BCD2010} provide standard $O(h^r)$ convergence estimates over the space of affine pyramids (recall the discussion at the end of Section~\ref{sec:bkgd}) for their respective elements.
Since the shape functions for $\calY_r\Lambda^0$ contain all the degree $r$ polynomials, as was just shown in Theorem~\ref{thm:polyn-rep}, an $O(h^r)$ estimate should hold for these elements over any mesh involving affine pyramids linked to $\calS_r\Lambda^0$ elements on affinely mapped hexahedra.
A formal study of such estimates is a topic for future work.

\section{Dimension optimality and conclusions}
\label{sec:dim-opt}

We compare the dimensions of $\calY_r\Lambda^0$ and $\calY_r^-\Lambda^0$ to other order $r$ pyramid elements in the literature, as summarized in Table~\ref{tab:dim-comp}.
In~\cite{NP2012}, Nigam and Phillips defined a set of shape functions, denoted $\calU^{(0),r}$, which has dimension $r^3+3r+1$.
Fuentes et al implemented the $\calU^{(0),r}$ shape functions, as described in~\cite{FKDN2015}.
In~\cite{BCD2010}, Bergot, Cohen, and Durufl\'e defined a space of rational functions that they denote $\hat P_r$, which has dimension $(r +1)(r +2)(2r +3)/6$, the same as $\dim \calY_r^-\Lambda^0$.
In~\cite{NP2012a}, Nigam and Phillips defined a reduced space, denoted $\calR^{(0)}_r$, which is a subset of $\calU^{(0),r}$, uses the $\calQ^{[r,r]}_r$ notation, and also has dimension equal to $\dim \calY_r^-\Lambda^0$.
Thus, there is likely very little practical distinction among $\calR^{(0)}_r$, $\hat P_r$, and  $\calY_r^-\Lambda^0$.

\begin{table}[h]
\begin{center}
\begin{tabular}{c|lllllll||l|l}
$r$                      & 1 & 2 & 3 & 4 & 5 & 6 & 7 & stated formula & reference\\[1mm]
\hline &&&&&&&&\\[-3mm]
$\dim\calY_r\Lambda^0$   & 5 & 13 & 25 & 42 & 65 & 95 & 133 & $(r^3 + 6 r^2 + 23 r)/6$\\[1mm]
$\dim\calY_r^-\Lambda^0$ & 5 & 14 & 30 & 55 & 91 & 140 & 204 & $(2 r^3 + 9 r^2 + 13 r + 6)/6$ \\[1mm]
$\dim\hat P_r$, $\dim\calR^{(0)}_r$ & 5 & 14 & 30 & 55 & 91 & 140 & 204 & (r +1)(r +2)(2r +3)/6 & \cite{BCD2010,NP2012a} \\[1mm]
$\dim\calU^{(0),r}$ &  5 & 15 & 37 & 77 & 141 & 235 & 365 & $r^3+ 3 r + 1$ & \cite{NP2012,FKDN2015}
\end{tabular}
\end{center}
\caption{Comparison of dimension counts for various pyramid elements in the literature.}
\label{tab:dim-comp}
\end{table}

\noindent
The space $\calY_r\Lambda^0$ is clearly distinct and of smaller dimension than any other order $r$ pyramid elements in the literature.
It cannot ``replace'' any of these elements, however, as it is has degrees of freedom that match $\calS_r\Lambda^0$ on a parallelogram face, not $\calQ_r^-\Lambda^0$.
Accordingly, the $\calY_r\Lambda^0$ element is not only useful but \textit{required} to create a conforming finite element method on hybrid meshes that employ serendipity elements on hexahedra.

In addition, the dimension of $\calY_r\Lambda^0$ agrees with the ``expected'' minimal dimension of a conforming finite element space on pyramids.
In~\cite{CG2015}, we discuss \textit{minimal compatible finite element systems} and state a criterion (Corollary 3.2) for computing the smallest possible dimension among all conforming finite element spaces that contain a desired set of functions, typically polynomials of degree at most $r$.
We show that the $\calP_r\Lambda^0~(=\calP_r^-\Lambda^0)$ spaces have this property on tetrahedra and the $\calS_r\Lambda^0$ spaces have this property on $n$-cubes\footnote{The claim about the minimality of $\calS_r\Lambda^0(\square_n)$ assumes that the remaining spaces in an exact sequence starting with $\calS_r\Lambda^0(\square_n)$ have decreasing polynomial approximation power.}.
Thus, it is expected that a minimal dimension element on a pyramid should correspond to $\calS_r\Lambda^0$ on its parallelogram face.

\begin{figure}
\begin{center}
\sbox{\strutbox}{\rule{0pt}{0pt}}           
\begin{tabular}[.8\textwidth]{@{\extracolsep{\fill}} cccc}
\parbox{.22\textwidth}{\includegraphics[width=.22\textwidth]{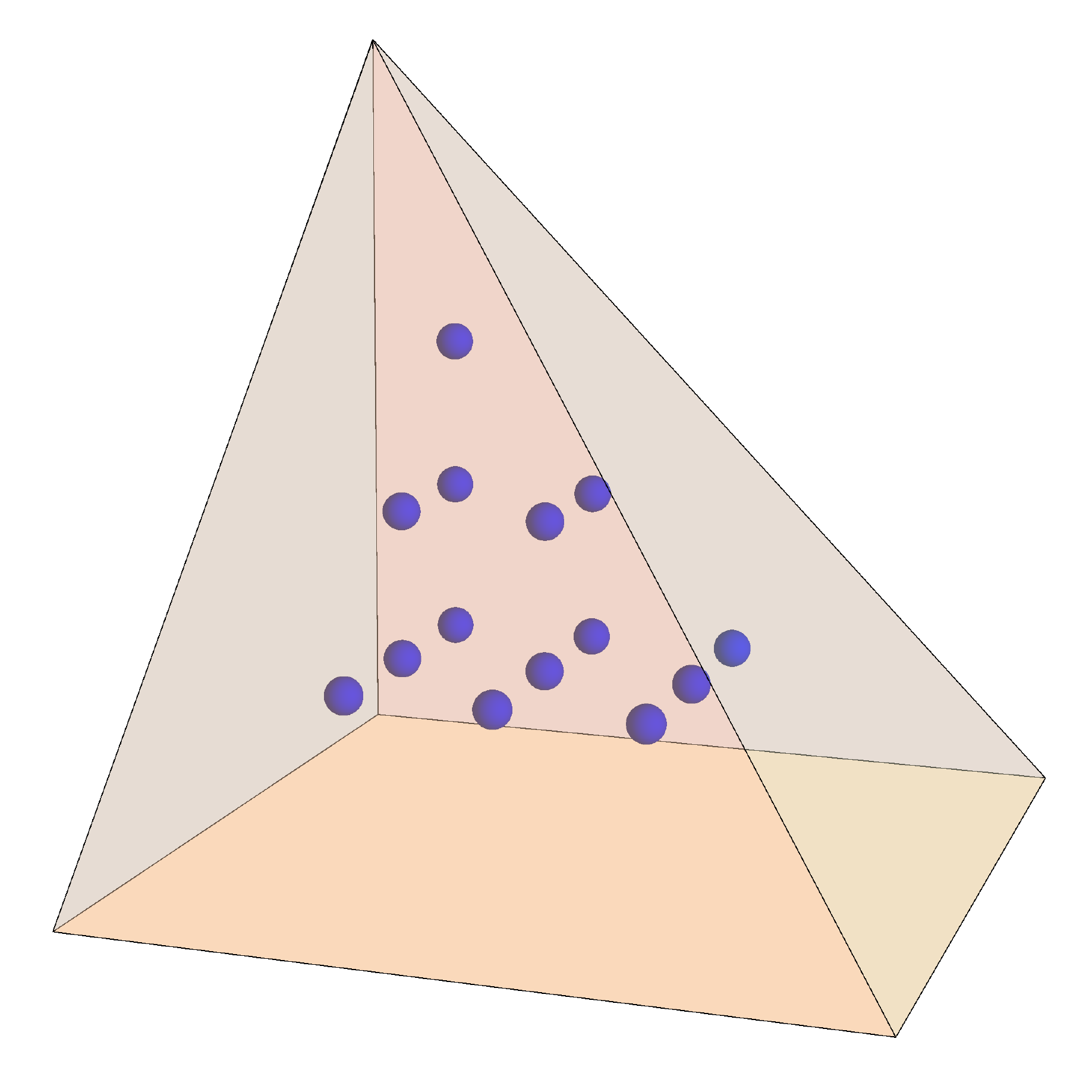}} & 
\parbox{.22\textwidth}{\includegraphics[width=.22\textwidth]{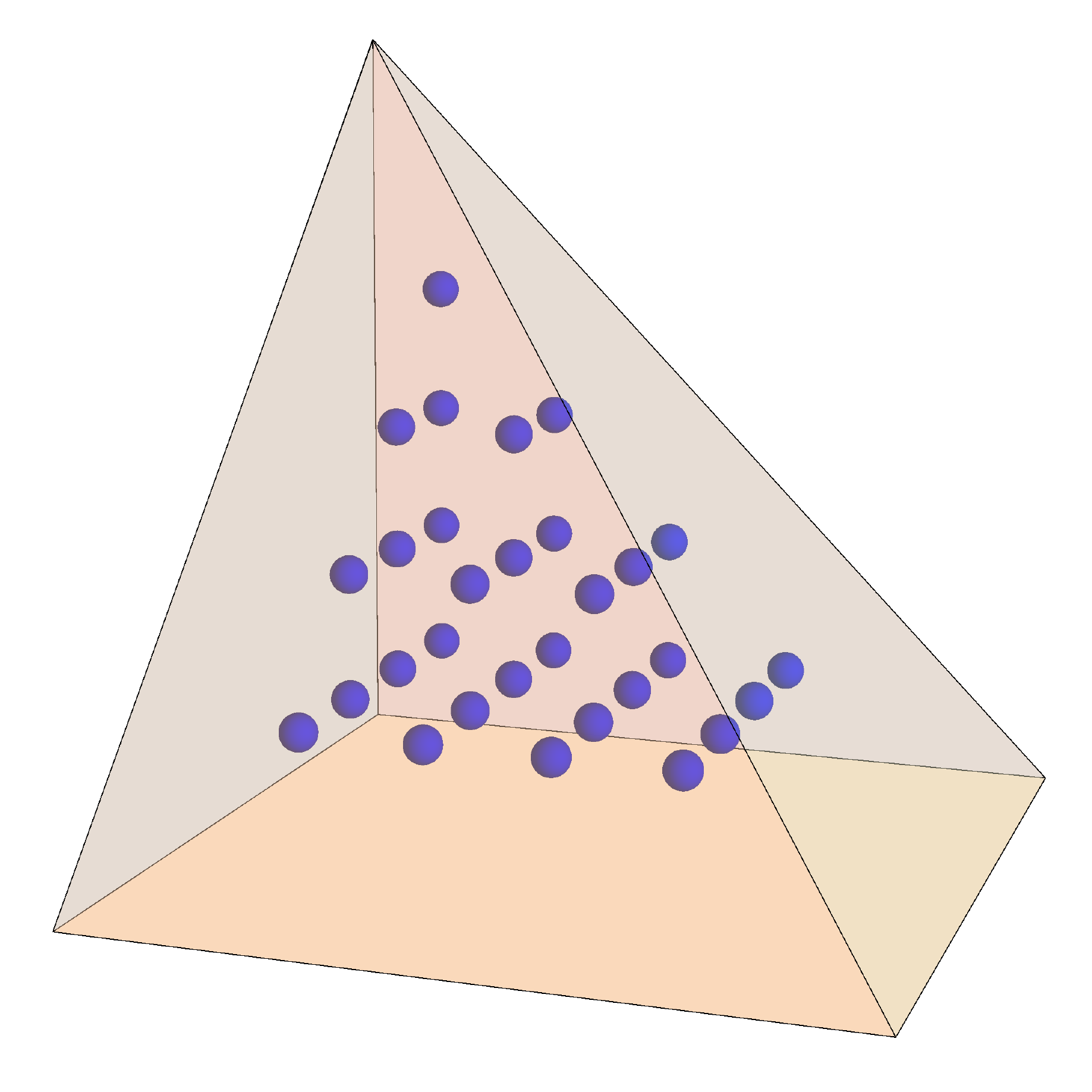}} &
\parbox{.22\textwidth}{\includegraphics[width=.22\textwidth]{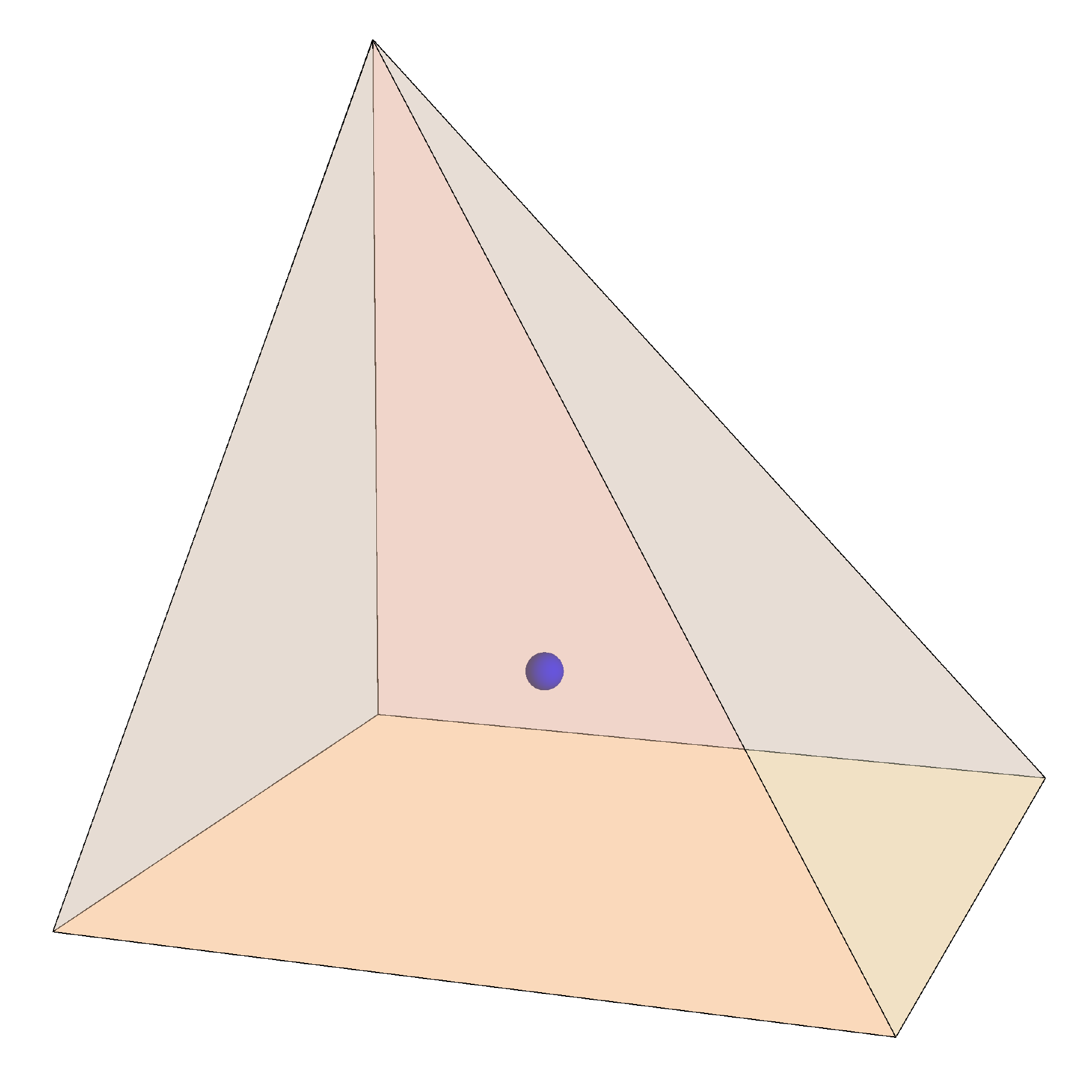}} & 
\parbox{.22\textwidth}{\includegraphics[width=.22\textwidth]{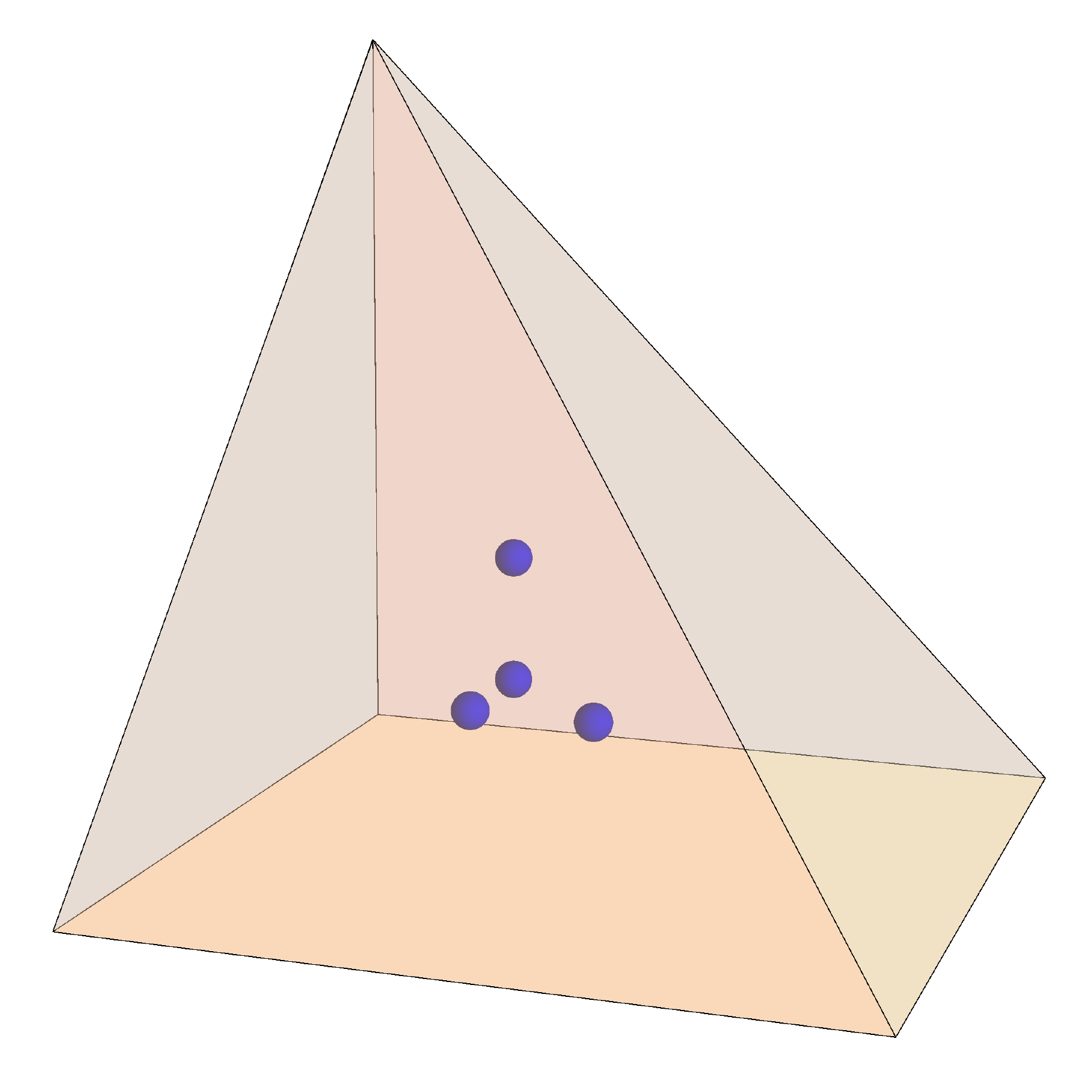}} 
\end{tabular}
\end{center}
\caption{
Degrees of freedom associated to the interior of the pyramid for $\calY_5^-\Lambda^0$, $\calY_6^-\Lambda^0$, $\calY_5\Lambda^0$, and $\calY_6\Lambda^0$ (left to right).
Note that there are no interior degrees of freedom for $\calY_r\Lambda^0$ when $r<5$.
}
\label{fig:int-dofs}
\end{figure}

The key question thus becomes the number of degrees of freedom that should be associated to the interior of $\hat K$.
For this, we note that the dimension of $R_\pyr$ for $\calY_r\Lambda^0$ is
\[\dim \phi\left(\bub\cdot\calS^{[r-5,r-5]}_{r-5}\right)\Lambda^3(\hat K) = {r-2 \choose 3} = \dim\calP_{r-5}\Lambda^3(\hat K).\]
We can map $\calP_{r-5}\Lambda^3(\hat K)$ bijectively to $\calP_{r}\Lambda_0^0(\hat K)$, the space of polynomials that vanish on the boundary of $\hat K$, by the map\\[-2mm]
\[p\,dV\longmapsto p\, \xi \eta \zeta (\xi+\zeta-1) (\eta+\zeta-1),\]
since $\xi \eta \zeta (\xi+\zeta-1) (\eta+\zeta-1)$ is a degree 5 polynomial that vanishes on $\p\hat K$.
Thus, $\calP_{r-5}\Lambda^3(\hat K)$ has a natural correspondence with the space of bubble functions of degree at most $r$ on pyramids.
Moreover, the growing field of virtual element methods, a distinct but related approach to the definition of finite element type methods on meshes with mixed mesh geometries, also assigns at least  ${r-2 \choose 3}$ degrees of freedom to the interior of a pyramid geometry~\cite{dVBMR2015}.
Therefore, it seems unlikely that the dimension of $\calY_r\Lambda^0$ could be reduced further without some loss to numerical accuracy or polynomial approximation order.

\paragraph{Conclusions and future work.}
In this work, we have used tools from finite element exterior calculus to define and analyze pyramidal finite elements that are affinely mapped from a reference geometry.
In particular, the use of superlinear degree in the description of serendipity elements was essential to the definition of the new serendipity-linking family, $\calY_r\Lambda^0$.
A number of topics related to the $\calY_r\Lambda^0$ family remain open for exploration, including the definition of local basis functions, efficient implementaion schemes, and integration with existing finite element solvers.
In light of the significant current interest in pyramid elements, progress in these areas is likely to occur very rapidly.

\text{}\\[0mm]
\textbf{Acknowledgements.}
This research was supported in part by NSF Award 1522289.
\vspace{-3mm}

\bibliographystyle{abbrv}
\bibliography{bib-pyr}

\appendix
\section{Additional background on finite element exterior calculus}
\label{appdx:feec}

Finite element exterior calculus uses tools from differential geometry and topology to define, classify, and analyze families of finite elements.
The scope of the theory is quite broad so we focus here only on those aspects that are immediately relevant to $H^1$-conforming finite element methods on pyramid geometries.

\paragraph{Polynomial differential $0$-forms and $n$-forms.}
Let $\Omega_n\subset\R^n$, $n\geq 1$, be an $n$-dimensional geometrical object in a finite element mesh.
The space $\calP_r\Lambda^0(\Omega_n)$ is defined to be the space of polynomials in $n$-variables of degree at most $r$.
The index `0' indicates that the elements of the space are differential 0-forms, i.e.\ scalar-valued functions.
Let $dV:=dx_1\cdots dx_n$ denote the volume form on $\Omega_n$, e.g.\ $dV= dxdy$ on a domain in $\R^2$ and $dV=dxdydz$ on a domain in $\R^3$.
The space $\calP_r\Lambda^n(\Omega_n)$ is then defined by
\[\calP_r\Lambda^n(\Omega_n) := \left\{~q dV~:~q\in\calP_r\Lambda^0(\Omega_n)~\right\}.\]
\indent The subtle but essential perspective of differential geometry is that a polynomial $p\in\calP_r\Lambda^0(\Omega_n)$ \textit{cannot} be integrated on $\Omega_n$ since it does not have the volume form attached.
On the other hand, a differential $n$-form $qdV\in\calP_r\Lambda^n(\Omega_n)$ can be integrated over $\Omega_n$, according to standard multivariate calculus techniques.
While seemingly pedantic, this perspective is ingrained even in first semester calculus where students can lose points for forgetting to write `$dx$' at the end of an integrand.
Spaces of differential $k$-forms for integers $k$ with $0<k<n$ require some more mathematical machinery to define, but are not needed in this work; definitions can be found in~\cite{AFW2006,AFW2010} or any textbook on differential geometry.

\paragraph{Trace operator.}
The trace operator associated to a subset $f\subset\Omega_n$ is a function $\tr_f:\calP_r\Lambda^0(\Omega_n)\raw\calP_r\Lambda^0(f)$.
The value of $\tr_f p$ is defined to be the pullback of $p$ via the inclusion map $f\hookrightarrow\Omega_n$, which can be interpreted as the restriction of $p$ to the domain $f$.
The trace operator is used to prove that a finite element family is $H^1$-conforming as follows.
If two elements $\Omega_n$ and $\Omega_n'$ meet along an $(n-1)$ dimensional face $f$ in a mesh, $H^1$-conformity requires that $\tr_f\calP_r\Lambda^0(\Omega_n)=\tr_f\calP_r\Lambda^0(\Omega_n')$.
For instance, let $\Omega_n=\triangle_2$ and $\Omega_n'=\triangle_2'$, where $\triangle_2$, $\triangle_2'$ are two triangles meeting along an edge $\be$ in a mesh.
Since restricting a polynomial from a plane to a line decreases its degree by 1, we have $\tr_f\calP_r\Lambda^0(\triangle_2)= \calP_{r-1}\Lambda^0(\be) =\tr_f\calP_r\Lambda^0(\triangle_2')$.
Additional details on trace operators in the context of pyramid finite elements can be found in~\cite{NP2012,NP2012a}.

\paragraph{Degrees of freedom.}
Given a space of shape functions on $\Omega_n$, degrees of freedom in a classical finite element sense are a set of functionals on $\Omega_n$ that take the shape functions as inputs.
Here, our shape functions are spaces of $0$-forms and our degrees of freedom require integration over $d$-dimensional portions of the pyramid geometry, for $d=0,1,2,3$.
Accordingly, each degree of freedom requires the input $u$ to be restricted to a differential $0$-form on the $d$-dimensional geometry portion, via an appropriate trace operator, and then multiplied by a differential $d$-form.
In general, a space of differential $k$-forms used as shape functions is indexed by a space of differential $d-k$ forms on a $d$-dimensional geometry and multiplication is generalized for $k>0$ via the wedge product.

\end{document}